\documentclass[12pt,a4paper]{article}
\usepackage{amsfonts}
\usepackage[pdftex]{graphicx}
\newtheorem{theorem}{Theorem}[section]
\newtheorem{lemma}{Lemma}[section]
\newtheorem{prop}{Proposition}[section]

\newcommand{\LP}{L\'{e}vy process}
\newcommand{\R}{\mathbb{R}}
\newcommand{\Rd}{\mathbb{R}^d}
\newcommand{\nN}{n \in \mathbb{N}}

\newcommand{\E}{\mathbb{E}}

\newcommand{\cadlag}{c\`adl\`ag}

\newcommand{\tr}{\mbox{tr}}

\newcommand{\bean}{\begin{eqnarray*}}
\newcommand{\eean}{\end{eqnarray*}}

\newcommand{\tra}{{\mathcal T}}

\date{}

\title{The Kalman-Bucy Filter for Integrable L\'{e}vy Processes with Infinite Second Moment}

\author{David Applebaum\footnote{D.Applebaum@sheffield.ac.uk}, Stefan Blackwood\footnote{stefan.blackwood@sheffield.ac.uk}\\
 School of Mathematics and Statistics,\\ University of Sheffield, \\
Sheffield S3 7RH\\
United Kingdom.}

\begin{document}
\maketitle

\begin{abstract} We extend the Kalman-Bucy filter to the case where both the system and observation processes are driven by finite dimensional L\'{e}vy processes, but whereas the process driving the system dynamics is square-integrable, that driving the observations is not; however it remains integrable.  The main result is that the components of the observation nose that have infinite variance make no contribution to the filtering equations. The key technique used is approximation by processes having bounded jumps.
\end{abstract}

\section{Introduction} The Kalman-Bucy filter is the continuous time version of the famous Kalman filter, which was developed in discrete time to find the ``best estimate'' of a linear system based on observations which are also assumed to dynamically evolve in a linear fashion. These techniques have had an enormous influence on control theory and signal processing, as well as the analysis of time series, and have found important applications within such areas as  navigational and guidance systems and satellite orbit determination (see e.g. \cite{Cip}.)

In this paper we focus on the continuous time Kalman-Bucy filter. This is usually presented with Brownian motion as driving noise (see e.g. \cite{Ox} chapter 6, \cite{BC} section 6.2 or \cite{Xi} Chapter 9), although it can be more generally developed when the noise has orthogonal increments \cite{Dav}, or is a fractional Brownian motion \cite{KB}. But all of these extensions require that both the signal and the observation noise have a finite second moment for all time. An exception to this is the paper \cite{BM}, but that uses a very different approach to the one that we will present, and it seems that the only concrete examples that fit naturally into that framework are the $\alpha$-stable L\'{e}vy processes. We also mention \cite{AF} that deals with a very similar problem to ours. We compare the approach of that paper with ours below. An extension of the discrete time Kalman filter to deal with some noise distributions having heavy tails may be found in \cite{SI}.

L\'{e}vy processes are the most general class of stochastic processes with suitably regular sample paths that have stationary and independent increments. Although they have been part of the classical theory of probability and stochastic processes since the 1930s, they are currently going through a stage of intense development and broad application due to their flexibility in modeling phenomena that combine continuous random motion with discontinuous jumps of arbitrary size. There have been many interesting applications of these processes, in e.g. mathematical finance (see Chapter 5 of \cite{App} and \cite{CoTa}), evolutionary biology \cite{Bar} and signal processing \cite{KP, App1}.

In this paper we describe how the Kalman-Bucy filter can be extended to L\'{e}vy noise in both the system and observations, with the system noise having a finite second moment while the observation noise has infinite second moment but finite mean. Our method is to truncate the size of the jumps in the observation noise, set up standard Kalman-Bucy filters for processes having finite second moments, and then take a passage to the limit. We are thus able to obtain a standard form of the Kalman-Bucy filter for a natural $L^{1}$ linear estimate of the system; note however that the solution of the  Riccati equation cannot any longer be interpreted as a ``mean square error'', although it is a limit of such terms. The components of the observation noise that have infinite variance make no contribution to either the Riccati equation or the Kalman-Bucy filter. In particular, if the variance of every component blows up to infinity (and this always happens in the one-dimensional case in this context) then the limiting Riccati equation linearises, the Kalman-Bucy filter is deterministic, and is the solution of a linear first-order differential equation.

 In contrast to our linear filter, \cite{AF} obtain a non-linear optimal filter for the {\it best measurable estimate}. Having the system evolve with Gaussian noise is also a vital aspect of their set-up, indeed the authors describe this as a ``significant limitation'' in the closing paragraph of their paper. We do not discuss optimality of the filter in this paper; this seems to be quite a difficult problem in our Banach space setting, where we cannot exploit orthogonality as in the standard Hilbert space framework.

The plan of the paper is as follows. In section 2 we gather together all the facts that we need about L\'{e}vy processes. In section 3, we set up the standard Kalman-Bucy filter for L\'{e}vy processes having finite second moments. In section 4 we describe the limiting processes, and obtain our main results. Some numerical simulations are presented in section 5, where we make practical comparisons of our filter with those of \cite{AF} and \cite{BM}.   A more detailed account of the results presented here can be found in \cite{Blac}. Finally we remark that non-linear filtering with L\'{e}vy processes has also recently received some development, and we direct the interested reader to \cite{MBP, PS, Blac}.

\vspace{5pt}

{\it Notation.} In this paper we will work  with vectors in $\R^{d}$ where $d$ takes at least four different values. If $v \in \R^{d}$ we denote its Euclidean norm by $|v|$ or $|v|_{d}$ if we want to emphasise the dimension of the space in which $v$ lies.  If $A$ is a $d \times m$ real valued matrix, its transpose will be denoted by $A^{\tra}$. For such matrices the operator norm is defined to be $||A||_{\infty} = \sup\{|Au|_{d}; |u|_{m} = 1\}$ and we frequently use the fact that if $u \in \R^{m}$ then $|Au|_{d} \leq ||A||_{\infty}|u|_{m}$. In the last part of the paper we will also need the Hilbert-Schmidt norm $||A||_{hs}:=\tr(A^{\tra}A)^{\frac{1}{2}}$ (where tr denotes the trace) and we recall the well known inequality $||A_{1}A_{2}||_{hs} \leq ||A_{1}||_{\infty}||A_{2}||_{hs}$. 
If  $(\Omega, {\cal F}, P)$ is a probability space, we use $||\cdot||_{p}$ to denote the norm in $L^{p}(\Omega, {\cal F}, P;\R^{d})$ for $1 < p < \infty$, so that if $X  = (X_{1}, \ldots, X_{d}) \in L^{p}(\Omega, {\cal F}, P;\R^{d})$, then $||X||_{p} = \E(|X|^{p})^{\frac{1}{p}}$. In this paper we will only be interested in the cases $p=1,2$ and we note that $L^{2}(\Omega, {\cal F}, P;\R^{d})$ is continuously embedded into $L^{1}(\Omega, {\cal F}, P;\R^{d})$ and that the embedding is a contraction. If $P$ is an orthogonal projection acting in a Hilbert space, we write $P^{\bot}:=I - P$, where $I$ is the identity operator.

\section{Preliminaries on L\'{e}vy Processes}

Let $(\Omega, {\cal F}, P)$ be a probability space and $({\cal F}_{t}, t \geq 0)$ be a filtration of ${\cal F}$ that satisfies the ``usual hypotheses'' of right continuity and completeness. Let $L = (L(t), t \geq 0)$ be a \LP~taking values in $\R^{d}$, that is defined on this space and adapted to the given filtration. This means that $L$ has stationary and independent increments in the strong sense that $L(t) - L(s)$ is independent of ${\cal F}_{s}$ for all $0 \leq s < t < \infty$, that $L$ is stochastically continuous and has \cadlag~paths, and that $L(0) = 0$~(a.s.). Under these conditions, we have the {\it L\'{e}vy-It\^{o} decomposition} (see \cite{App}, section 2.4, pp.112-31)
\begin{equation} \label{LI}
L(t) = bt + W_{a}(t) + \int_{|y| \leq 1}y\tilde{N}(t, dy) + \int_{|y| > 1}yN(t,dy)
\end{equation}
for all $t \geq 0$. Here $b \in \Rd, W_{a}$ is a Brownian motion in $\Rd$ with Cov$(W^{i}_{a}(s), W^{j}_{a}(t)) = (s \wedge t)a_{ij}$ for all $1 \leq i,j \leq d, s, t \geq 0$ where $a = (a_{ij})$ is a non-negative definite symmetric $d \times d$ matrix, $N$ is a Poisson random measure defined on $\R^{+} \times \R^{d}$ which is independent of $W_{a}$, and $\tilde{N}$ is its compensator for which $\tilde{N}(dt,dy) = N(dt,dy) -dt \nu(dy)$, where $\nu$ is a L\'{e}vy measure on $\Rd$, i.e. a Borel measure for which $\nu(\{0\}) = 0$ and $\int_{\Rd}(|y|^{2} \wedge 1) \nu(dy) < \infty$. Note that we can always write $W_{a}^{i}(t) = \sum_{j=1}^{m}\sigma^{i}_{j}W^{j}(t)$ where $W = (W^{1}, \ldots, W^{m})$ is a standard Brownian motion in $\R^{m}$ and $\sigma = (\sigma^{i}_{j})$ is a $d \times m$ matrix for which $\sigma \sigma^{\tra} = a$.

We say that a \LP~$L$ is $p$-integrable for $p > 0$ if $\E(|L(t)|^{p}) < \infty$ for all $t \geq 0$. It is well-known (see e.g. \cite{App}, section 2.5, pp.131-3) that $L$ is $p$-integrable for $p \in \mathbb{N}$ if and only if $\int_{|y| > 1}|y|^{p}\nu(dy) < \infty$.  In particular we will frequently deal with integrable L\'{e}vy processes ($p=1$) which are also centred, i.e. $\E(L(t)) = 0$ for all $t \geq 0$. In this case the L\'{e}vy-It\^{o} decomposition can be written
\begin{equation} \label{LI1}
L(t) = W_{a}(t) + \int_{\R^{d}}y\tilde{N}(t, dy).
\end{equation}
If such a process $L$ is also square-integrable ($p=2$) we have Cov$(L^{i}(s), L^{j}(t)) = (s \wedge t)\Theta_{ij}$ for all $1 \leq i,j \leq d, s, t \geq 0$, where $\Theta_{ij} = a_{ij} + \psi_{ij}$ and
\begin{equation} \label{Cov}
\psi_{ij} = \int_{\R^{d}}y^{i}y^{j}\nu(dy).
\end{equation}
We always assume that the matrix $\psi = (\psi_{ij})$ is non-degenerate and so has strictly positive determinant.
We will also need $\lambda = \tr(\Theta)$ and note that
\begin{equation} \label{Cov1}
\lambda = \tr(a) + \int_{\R^{d}}|y|^{2}\nu(dy).
\end{equation}

Next we describe an approximation result that we will find useful in section 4. Let $L$ be a centred integrable L\'{e}vy process as in (\ref{LI1}) and let $(L_{n}, \nN)$ be the sequence of centred square-integrable L\'{e}vy processes where each $L_{n}$ has the L\'{e}vy-It\^{o} decomposition (for $t \geq 0$)
\begin{equation} \label{LI2}
L_{n}(t) = W_{a}(t) + \int_{|y| \leq n}y\tilde{N}(t, dy).
\end{equation}

\begin{prop} \label{approx1} The sequence $(L_{n}(t), \nN)$ converges to $L(t)$ for all $t \geq 0$ in the $L^{1}$-sense and the convergence is uniform on finite intervals.
\end{prop}

{\it Proof.} For all $T \geq 0$, \bean
\sup_{0 \leq t \leq T}||L(t) - L_{n}(t)||_{1} & = &  \sup_{0 \leq t \leq T}\E\left(\left|\int_{|y| > n}y\tilde{N}(t, dy)\right|\right)\\
& \leq & 2T \int_{|y| \geq n}|y|\nu(dy) \rightarrow 0~\mbox{as}~n \rightarrow \infty, \eean
where we have used the well-known fact that if $h \in L^{1}(\R^{d}, \nu)$ then $$\E\left(\left|\int_{\R^{d}}h(y)\tilde{N}(t, dy)\right|\right) \leq 2t\int_{\Rd}|h(y)|\nu(dy).~~~~~~\hfill \Box$$

\vspace{5pt}

Now suppose that $L$ is integrable but not square-integrable. Then $\int_{|y| \geq 1}|y|^{2}\nu(dy) = \infty$ and so
$\int_{|y| \geq 1}y_{i}^{2}\nu(dy) = \infty$ for $q$ values of $i$ where $1 \leq i \leq d$. Assume for simplicity that these $q$ values are $1, 2, \ldots, q$ and write  $\beta_{\nu}(n): = \max\left\{\int_{|y| \leq n}y_{i}^{2}\nu(dy); 1 \leq i \leq q\right\}$. In order to develop a manageable theory, we will also make the assumption that $\left|\int_{|y| \geq 1}y_{i}y_{j}\nu(dy)\right| < \infty$ for all $i \neq j$. An example of such a set up arises when the components of $L$ are independent, centred one-dimensional L\'{e}vy processes, $q$ of which have infinite variance. In that case, for each $i \neq j$,
$$ \int_{|y| \geq 1}y_{i}y_{j}\nu(dy) = \left(\int_{|y| \geq 1}y_{i}\nu(dy)\right)\left(\int_{|y| \geq 1}y_{j}\nu(dy)\right) = 0.$$
More generally, we say that $L$ is of {\it fully infinite variance} whenever $q = d$.

For each $\nN$, let $\Lambda^{(n)} = (\Lambda^{(n)}_{ij})$ be the covariance matrix corresponding to the jump part of the process $L^{n}$, as in (\ref{Cov}). In section 4, we will require the existence of $\lim_{n \rightarrow \infty}(\Lambda^{(n)})^{-1}$. If $d = 1$ it is easy to see that this limit exists and is zero. If $d > 1$ and $q=d$, then each element of the matrix of cofactors of $\Lambda^{(n)}$ is $O(\beta_{\nu}(n)^{\alpha})$ as $n \rightarrow \infty$, where $1 \leq \alpha \leq d-1$, however $\det(\Lambda^{(n)})$ is $O(\beta_{\nu}(n)^{d})$. It follows that $\lim_{n \rightarrow \infty}(\Lambda^{(n)})^{-1} = 0$ also in this case.

If $d > 1$ and $q < d$ then the sequence $(\Lambda^{(n)}{^{-1}}, \nN)$ may converge to a non-zero limit as $n \rightarrow \infty$. For an illustrative example, let $(e_{i}, 1 \leq i \leq d)$ be the natural basis in $\R^{d}$ and for each $\nN$, let $L_{n}$ be the L\'{e}vy process:
$$ L_{n}(t) = \sum_{i=1}^{q}J_{n}^{(i)}(t)e_{i} + \sum_{i=q+1}^{d}W^{(i)}(t)e_{i},$$
where $t \geq 0, W^{(q+1)}, \ldots, W^{(d)}$ are independent one-dimensional standard Brownian motions and $J_{n}^{(1)}, \ldots, J_{n}^{(q)}$ are independent one-dimensional pure jump L\'{e}vy processes, which are also independent of all the $W^{(i)}$'s, with each
$ J_{n}^{(i)}(t) = \int_{|y_{i}| < n}y\tilde{N_{i}}(t, dy_{i})$ having L\'{e}vy measure that is concentrated on the linear span of the vector $e_{i}$.
In this case, for each $\nN$,
$$ \Lambda^{(n)} = \mbox{diag}\left(\int_{y_{1}^{2} < n}y_{1}^{2}\nu_{1}(dy), \ldots, \int_{y_{q}^{2} < n}y_{q}^{2}\nu_{q}(dy), 1, \ldots , 1\right),$$
and we have
$$ \lim_{n \rightarrow \infty}(\Lambda^{(n)})^{-1} = \mbox{diag}\left(0, \ldots,0, 1, \ldots , 1\right).$$

In the sequel, we will always assume that $ \lim_{n \rightarrow \infty}(\Lambda^{(n)})^{-1}$ exists, and we denote the limit by $\Upsilon_{\infty}$.

\section{The Kalman-Bucy Filter With Square-Integrable L\'{e}vy Processes}

The Kalman-Bucy filter has been extensively developed in the case where the noise (at both the system and the observation level) is in the general class of processes with orthogonal increments (see e.g. \cite{Dav}) and the theory we describe in this section is no more than a special case of this. However it is important for us to write down the key equations as we will need to use these in the main part of the paper where we go beyond the framework of square-integrability. We also point out that in the literature (see also e.g. \cite{Ox, BC, Xi}), it is typical for authors to present the theory using standard Brownian motion, or a noise process whose covariance matrix is the identity. When we use L\'{e}vy processes, the covariance matrix is typically non-trivial (see (\ref{Cov})) and it will play a key role in the approximation scheme that we develop in the next section.

We begin by formulating the general linear filtering problem within the context of L\'{e}vy processes. Let $L_{1}$ and $L_{2}$ be independent square-integrable L\'{e}vy processes defined on $\Omega$ and taking values in $\R^{l}$ and $\R^{p}$ (respectively). In the sequel we will use the notation $b_{i}, \nu_{i}, \lambda_{i}$ etc. when dealing with characteristics of L\'{e}vy processes $L_{i} (i = 1,2)$ as described in section 2.

Consider the following stochastic differential equation (SDE) in $\R^{d_{1}}$:
\begin{equation} \label{sys}
dY(t) = A(t)Y(t-)dt + B(t)dL_{1}(t),
\end{equation}
Here $A$ and $B$ are locally bounded left continuous functions taking values in the space of $d_{1} \times d_{1}$ and $d_{1} \times l$ matrices (respectively). We assume that the initial condition $Y_{0}: = Y(0)$ is a ${\cal F}_{0}$-measurable random variable. Under these conditions, the SDE (\ref{sys}) has a unique \cadlag~solution (see e.g. \cite{App}, Chapter 6) $Y = (Y(t), t \geq 0)$ taking values in $\R^{d_{1}}$. We call $Y$ the {\it system process}. Note that using the variation of constants formula we have
\begin{equation} \label{sys1}
Y(t) = \exp{\left\{\int_{0}^{t}A(u)du\right\}}Y(0) + \int_{0}^{t}\exp{\left\{\int_{s}^{t}A(u)du\right\}}B(s)dL_{1}(s),
\end{equation}
for all $t \geq 0$.  We have that $\E(|Y(t)|^{2}) < \infty$ for all $t \geq 0$ (see e.g. Corollary 6.2.4 on p.373 of \cite{App}.) Now let $d_{2} \leq d_{1}$ and let $C$ and $D$ be locally bounded, left continuous functions taking values in the space of $d_{2} \times d_{1}$ and $d_{2} \times p$ matrices (respectively). Then the SDE
\begin{equation}  \label{obs}
dZ(t) = C(t)Y(t)dt + D(t)dL_{2}(t)
\end{equation}
with initial condition $Z(0) = 0$ (a.s.) has a unique \cadlag~solution $Z = (Z(t), t \geq 0)$  taking values in $\R^{d_{2}}$. We call $Z$ the {\it observation process}. Just as in the case of the system, $\E(|Z(t)|^{2}) < \infty$ for all $t \geq 0$.

From now on we will restrict ourselves to working on a time interval $[0,T]$. We will also assume that the $d_{2} \times d_{2}$ matrix $D(t)D(t)^{\mathcal T}$ is invertible for all $0 \leq t \leq T$ and that the mapping $t \rightarrow D(t)D(t)^{\mathcal T}$ is bounded away from zero. We define $G(t):= (D(t)D(t)^{\mathcal T})^{-\frac{1}{2}}$ for all $0 \leq t \leq T$.

Let $X = (X(t), 0 \leq t \leq T)$ be a square-integrable stochastic process defined on $\Omega$ and taking values in $\R^{d_{2}}$. Let ${\cal L}(X)$ denote the closure in $L^{2}(\Omega, {\cal F},P; \R^{d_{1}})$ of all finite linear combinations of random variables of the form $$c_{0} + c_{1}X(t_{1}) + \cdots + c_{n}X(t_{n})$$ where $c_{0} \in \R^{d_{1}}, c_{1}, \ldots, c_{n}$ are arbitrary $d_{1} \times d_{2}$ matrices and $0 \leq  t_{1} < \cdots < t_{n} \leq T$. We will use $P_{X}$ to denote the orthogonal projection in $L^{2}(\Omega, {\cal F},P; \R^{d_{1}})$ whose range is ${\cal L}(X)$. Now for all $0 \leq t \leq T$, we define
$$ \widehat{Y}(t) = P_{Z}(Y(t)),$$
then $\widehat{Y}$ is the {\it best linear estimator} of $Y$. Note that $\E(\widehat{Y}(t)) = \E(Y(t))$ for all $0 \leq t \leq T$ and $\widehat{Y}(0) = \mu_{0}:= \E(Y_{0})$. We define the {\it innovations process} $N =(N(t), 0 \leq t \leq T)$ by the prescription:
$$ N(t) = Z(t) - \int_{0}^{t}C(s)\widehat{Y}(s-)ds,$$
then $N$ is a process with orthogonal increments taking values in $\R^{d_{2}}$. We further introduce the process $R = (R(t), 0 \leq t \leq T)$ taking values in $\R^{d_{2}}$ which is the unique solution of the SDE:
$$ dR(t) = G(t)dN(t),$$
with $R(0) = 0$ (a.s.). Then $R$ is centred and has orthogonal increments with covariance $\E(R(s)R(t)^{\mathcal T}) = \Sigma_{2}(s \wedge t)$ for all $0 \leq s, t \leq T$ where for $0 \leq u \leq T$,
$$ \Sigma_{2}(u): = \int_{0}^{u}G(r)D(r)\Lambda_{2}D(r)^{\mathcal T}G(r)^{\mathcal T}dr.$$
We assume that $\Sigma_{2}(t)^{-1}$ exists for all $0 \leq t \leq T$ and that the mapping $t \rightarrow \Sigma_{2}(t)^{-1}$ is bounded.
Note that we also have that $\E(R(s)^{T}R(t)) = \lambda_{2}(s \wedge t)$. We then find that for all $0 \leq t \leq T$,
we have the following representation of the estimator.
$$ \widehat{Y}(t) = \E(Y(t)) + \int_{0}^{t}\frac{\partial}{\partial s}\E(Y(s)R(t)^{\mathcal T})\Sigma_{2}(s)^{-1}ds.$$

We define the {\it mean squared error}
$$ S(t) = \E((Y(t) - \widehat{Y}(t))(Y(t) - \widehat{Y}(t))^{\mathcal T}),$$
for each $0 \leq t \leq T$. Then $S$ satisfies the {\it Riccati equation}
\begin{eqnarray} \label{Ricc}
\frac{dS(t)}{dt} & = & A(t)S(t) + S(t)A(t)^{\tra} + B(t)\Lambda_{1}B(t)^{\tra} \nonumber \\
& - & S(t)C(t)^{\tra}G(t)^{\tra}\Sigma_{2}(t)^{-1}G(t)C(t)S(t)^{\tra},
\end{eqnarray}
with initial condition $S(0) = \mbox{Cov}(Y_{0})$, and the {\it Kalman-Bucy filter} represents the best linear estimate $\widehat{Y}$ as the solution of the following SDE:
\begin{eqnarray} \label{KBFstand}
d\widehat{Y}(t) & = & A(t)\widehat{Y}(t)dt \nonumber \\
& + & S(t)C(t)^{\tra}G(t)^{\tra}\Sigma_{2}(t)^{-1}G(t)(dZ(t) - C(t)\widehat{Y}(t)dt).
\end{eqnarray}
Full details of the above derivations can be found in \cite{Blac}, but we stress that these results are essentially known (see e.g. \cite{Dav}).

\section{The Kalman-Bucy Filter With Integrable (Infinite Second Moment) L\'{e}vy Processes}

From now on we assume that the L\'{e}vy process $L_{2}$ is centred and integrable so that it has L\'{e}vy-It\^{o} decomposition (\ref{LI1}). We further assume that it fails to be square-integrable, i.e. $\E(|L_{2}(t)|^{2}) = \infty$ for all $t > 0$ and so $\int_{|y| > 1}|y|^{2}\nu(dy) = \infty$\footnote{From now on, we drop the subscript $_{2}$ for all characteristics and component processes associated with $L_{2}$.}. For each $\nN$ we consider the approximating L\'{e}vy processes $L_{2}^{(n)}$ that have L\'{e}vy-It\^{o} decomposition (\ref{LI2}). We now set up a sequence of observation processes $(Z_{n}, \nN)$ where for each $\nN, 0 \leq t \leq T$,
\begin{equation}  \label{obs1}
dZ_{n}(t) = C(t)Y(t-)dt + D(t)dL_{2}^{(n)}(t)
\end{equation}


\begin{prop} \label{approx2} The sequence $(Z_{n}(t), \nN)$ converges to $Z(t)$ for all $t \geq 0$ in the $L^{1}$-sense and the convergence is uniform on finite intervals.
\end{prop}

{\it Proof.} This is by a similar argument to that of Proposition \ref{approx1}. $\hfill \Box$

\vspace{5pt}

For simplicity we now write $L^{p}(\Omega):=L^{p}(\Omega, {\cal F},P; \R^{d_{1}}), p=1,2$. For all $\nN$, we can form the linear spaces ${\cal L}(Z_{n})$ as described in section 3. The orthogonal projection from $L^{2}(\Omega)$ to ${\cal L}(Z_{n})$ is denoted $P_{n}$ and for each $0 \leq t \leq T$, we write $\widehat{Y}_{n}(t) = P_{n}(Y(t))$. The mean square error of $\widehat{Y}_{n}(t)$ is denoted by $S_{n}(t)$.

\begin{lemma} \label{est}
For all $\nN, 0 \leq t \leq T$,
\begin{equation} \label{est1}
   ||S_{n}(t)||_{\infty} \leq d_{1}||Y(t)||^{2}_{2}
\end{equation}
\end{lemma}

{\it Proof.} For each $1 \leq  i,j \leq d_{1}$,
\bean  |S_{n}(t)_{ij}| & \leq & \E(|Y(t)_{i} - \widehat{Y}_{n}(t)_{i}||Y(t)_{j} - \widehat{Y}_{n}(t)_{j}|) \\
& \leq & \E(|Y(t) - \widehat{Y}_{n}(t)|^{2}_{d_{1}})\\
& = & ||P_{n}^{\bot}Y(t)||_{2}^{2}\\
& \leq & ||Y(t)||^{2}_{2}, \eean
and the result follows by a easy matrix estimate. $\hfill \Box$

\vspace{5pt}



 In order to proceed further, we assume that the sequence of matrices $(\Sigma_{2}^{(n)}(t)^{-1}, \nN)$ converges for each $t \in [0,T]$ to a matrix we denote by $\Phi(t)$. 
 For example, if for all $0 \leq t \leq T, D(t) = Df(t)$ where $D$ is invertible and $f$ is strictly positive, then $G(t) = Gf(t)^{-1}$ where $G = (DD^{\tra})^{-1}$, then
 $$ \Sigma_{2}^{(n)}(t) = tGD\Lambda_{2}^{(n)}D^{\tra}G^{\tra},$$

 and

 $$ \Phi(t) = \frac{1}{t}(G^{\tra})^{-1}(D^{\tra})^{-1}\Upsilon_{\infty} D^{-1}G^{-1}.$$ If the process $Z_{2}$ has fully infinite variance, then it is clear that $\Phi(t) = 0$, for all $0 \leq t \leq T$.

 In fact, we will need to impose a somewhat stronger condition, namely that $(G(t)^{\tra}\Sigma_{2}^{(n)}(t)^{-1}G(t), \nN)$ converges uniformly to $\Phi(t)$ on $(0, T]$. In the special case just described, this condition is satisfied if $f$ is such that $\inf_{t > 0} tf(t)^{2} > 0$. The search for general conditions for convergence seems to be a tricky problem in matrix analysis, and we will will not pursue it further here.

To simplify the presentation of proofs of the next two theorems, we introduce the notation $M(t): = G(t)C(t)$ for $0 \leq t \leq T$.

\begin{theorem} \label{Riccinf}
For each $0 \leq t \leq T$, the sequence $(S_{n}(t), \nN)$ converges to a matrix $S_{\infty}(t)$. The mapping $t \rightarrow S_{\infty}(t)$ is differentiable on $(0, T]$ and $S_{\infty}$ is the unique solution of the Riccati equation
\begin{eqnarray} \label{Riccinf1}
\frac{dS_{\infty}(t)}{dt} & = & A(t)S_{\infty}(t) + S_{\infty}(t)A(t)^{\tra} + B(t)\Lambda_{1}B(t)^{\tra} \nonumber \\
& - & S_{\infty}(t)C(t)^{\tra}G(t)^{\tra}\Phi(t)G(t)C(t)S_{\infty}(t)^{\tra},
\end{eqnarray}
\end{theorem}

{\it Proof.} Let $(\Xi(t); 0 \leq t \leq T)$ be the unique solution of the linear differential equation
\bean
\frac{d\Xi(t)}{dt} &  = &  A(t)\Xi(t) + \Xi(t)A(t)^{\tra} + B(t)\Lambda_{1}B(t)^{\tra}\\
& - &  \Xi(t)M(t)^{\tra}\Phi(t)M(t)\Xi(t)^{\tra},
\eean

with initial condition, $\Xi(0) = \mbox{Cov}(Y_{0})$. We will show that $\Xi(t) = S_{\infty}(t)$ for all $0 \leq t \leq T$.
Using (\ref{Ricc}), for all $\nN$ we have
\bean \frac{d}{dt}(\Xi(t)- S_{n}(t)) & = & A(t)(\Xi(t) - S_{n}(t)) + (\Xi(t)- S_{n}(t))A(t)^{\tra} \\
& + & S_{n}(t)M(t)^{\tra}\Sigma_{2}^{(n)}(t)^{-1}M(t)S_{n}(t)^{\tra}\\
& - &  \Xi(t)M(t)^{\tra}\Phi(t)M(t)\Xi(t)^{\tra} \\
& = & A(t)(\Xi(t) - S_{n}(t)) + (\Xi(t)- S_{n}(t))A(t)^{\tra} \\
& + & (S_{n}(t)-\Xi(t))M(t)^{\tra}\Sigma_{2}^{(n)}(t)^{-1}M(t)S_{n}(t)^{\tra}\\
& + & \Xi(t)M(t)^{\tra}(\Sigma_{2}^{(n)}(t)^{-1} - \Phi(t))M(t)S_{n}(t)^{\tra}\\
& + &  \Xi(t)M(t)^{\tra}\Phi(t)M(t)(S_{n}(t)^{\tra}-\Xi(t)^{\tra}).
\eean

Integrating both sides, taking matrix $||\cdot||_{\infty}$ norms and using (\ref{est}) yields for all $0 < t \leq T$,
\bean ||\Xi(t)- S_{n}(t)||_{\infty} & \leq & 2\int_{0}^{t}||A(r)||_{\infty}||\Xi(r)- S_{n}(r)||_{\infty}dr \\ & + & \int_{0}^{t}||S_{n}(r) - \Xi(r)||_{\infty}||C(r)||_{\infty}^{2}(||S_{n}(r)||_{\infty}||G(r)^{\tra}\Sigma_{2}^{(n)}(r)^{-1}G(r)||_{\infty}\\
& + & ||\Phi(r)||_{\infty}||\Xi(r)||_{\infty})dr\\
& + & \int_{0}^{t}||\Sigma_{2}^{(n)}(t)^{-1} - \Phi(t)||_{\infty}||M(r)||_{\infty}^{2}||\Xi(r)||_{\infty}||S_{n}(r)||_{\infty}dr.
\eean

By our assumptions, we have
$$\sup_{\nN}\sup_{0 < t \leq T}||G(t)^{\tra}\Sigma_{2}^{(n)}(t)^{-1}G(t)||_{\infty}^{2} = K <  \infty. $$


We now apply Gronwall's inequality and use (\ref{est1}) to obtain

\bean
||\Xi(t)- S_{n}(t)||_{\infty} & \leq & d_{1}\exp{\left\{\int_{0}^{t}[2||A(r)||_{\infty} + ||M(r)||_{\infty}^{2}(K||Y(r)||^{2}_{2} + ||\Phi(r)||_{\infty}||\Xi(r)||_{\infty})]dr\right\}}\\
& \times & \int_{0}^{t}||\Sigma_{2}^{(n)}(r)^{-1} - \Phi(r)||_{\infty}||M(r)||_{\infty}^{2}||\Xi(r)||_{\infty}||Y(r)||_{2}^{2}dr.
\eean

and so by dominated convergence, $\lim_{n \rightarrow \infty}||\Xi(t)- S_{n}(t)||_{\infty} = 0$, and the result follows.~~~~~~~~~$\hfill \Box$

\vspace{5pt}

Note that when the observation process $Z$ has fully infinite variance, then $\Phi(t)$ is identically zero for all $0 \leq t \leq T$ and (\ref{Riccinf1}) linearises and takes the form:

\begin{equation} \label{Riccinf2}
\frac{dS_{\infty}(t)}{dt}  = A(t)S_{\infty}(t) + S_{\infty}(t)A(t)^{\tra} + B(t)\Lambda_{1}B(t)^{\tra}.
\end{equation}

In particular, (\ref{Riccinf2}) always holds in the one-dimensional case.

\vspace{5pt}

Our next result is the desired $L^{1}$- L\'{e}vy Kalman-Bucy filter:

\begin{theorem} For each $0 \leq t \leq T$, the sequence $(\widehat{Y}_{n}(t), \nN)$ converges in $L^{1}(\Omega)$ to a random variable $\widehat{Y}(t)$.  The process $\widehat{Y}$ is the solution of the following SDE:
\begin{eqnarray} \label{KBFinf}
d\widehat{Y}(t) & = & A(t)\widehat{Y}(t)dt \nonumber \\
& + & S_{\infty}(t)C(t)^{\tra}G(t)^{\tra}\Phi(t)G(t)(dZ(t) - C(t)\widehat{Y}(t)dt)
\end{eqnarray}
\end{theorem}

{\it Proof.} Let $\Psi = (\Psi(t), 0 \leq t \leq T)$ be the unique solution (see e.g. Chapter 6 of \cite{App}) of the SDE
\bean d\Psi(t) & = & A(t)\Psi(t)dt \nonumber \\
& + & S_{\infty}(t)C(t)^{\tra}G(t)^{\tra}\Phi(t)G(t)(dZ(t) - C(t)\Psi(t)dt) \eean
with initial condition $\Psi(0) = \mu_{0}$ (a.s.) We will show that $\lim_{n \rightarrow \infty}||\Psi(t) - \widehat{Y_{n}}(t)||_{1} = 0$ and the result then follows by uniqueness of limits. We use a similar argument to that in the proof of Theorem \ref{Riccinf}. Using (\ref{KBFstand}) we find that
\bean \Psi(t) - \widehat{Y_{n}}(t) & = & \int_{0}^{t}A(r)(\Psi(r) - \widehat{Y_{n}}(r))dr\\
& + & \int_{0}^{t}S_{\infty}(r)M(r)^{\tra}\Phi(r)G(r)dZ(r)\\
& - & \int_{0}^{t}S_{n}(r)M(r)^{\tra}\Sigma_{n}(r)^{-1}G(r)dZ_{n}(r)\\
& - & \int_{0}^{t}S_{\infty}(r)M(r)^{\tra}\Phi(r)M(r)\Psi(r)dr\\
& + &  \int_{0}^{t}S_{n}(r)M(r)^{\tra}\Sigma_{n}(r)^{-1}M(r)\widehat{Y_{n}}(r)dr\\
& = & \int_{0}^{t}A(r)(\Psi(r) - \widehat{Y_{n}}(r))dr\\
& + & \int_{0}^{t}(S_{\infty}(r) - S_{n}(r))M(r)^{\tra}\Phi(r)G(r)dZ(r)\\
& + & \int_{0}^{t}S_{n}(r)M(r)^{\tra}(\Phi(r) - \Sigma_{n}(r)^{-1})G(r)dZ(r)\\
& + & \int_{0}^{t}S_{n}(r)M(r)^{\tra}\Sigma_{n}(r)^{-1}G(r)(dZ(r) - dZ_{n}(r))\\
& - & \int_{0}^{t}(S_{\infty}(r) - S_{n}(r))M(r)^{\tra}\Phi(r)M(r)\Psi(r)dr\\
& - &  \int_{0}^{t}S_{n}(r)M(r)^{\tra}(\Phi(r) - \Sigma_{n}(r)^{-1})M(r)\Psi(r)dr\\
& - & \int_{0}^{t}S_{n}(r)M(r)^{\tra}\Sigma_{n}(r)^{-1}M(r)(\Psi(r) - \widehat{Y_{n}}(r))dr.\eean
Taking the $L^{1}$-norm and using the fact that $Z(t) - Z_{n}(t)= \int_{0}^{t}\int_{|y| > n} D(r)y\tilde{N}(dr, dy)$ we find that
\begin{eqnarray} \label{last} ||\Psi(t) - \widehat{Y_{n}}(t)||_{1} & \leq & \int_{0}^{t}(||A(r)||_{\infty} + || S_{\infty}(r)M(r)^{\tra}\Sigma_{n}(r)^{-1}M(r)||_{\infty})||\Psi(r) - \widehat{Y_{n}}(r))||_{1}dr \nonumber \\
& + &  \E\left(\left|\int_{0}^{t}(S_{\infty}(r) - S_{n}(r))M(r)^{\tra}\Phi(r)G(r)dZ(r)\right|\right)\nonumber \\
& +  & \E\left(\left|\int_{0}^{t}S_{n}(r)M(r)^{\tra}(\Phi(r) - \Sigma_{n}(r)^{-1})G(r)dZ(r)\right|\right)\nonumber \\
& + & 2\int_{|y| > n}|y|\nu(dy)\int_{0}^{t}||S_{n}(r)M(r)^{\tra}\Sigma_{n}(r)^{-1})G(r)||_{\infty}dr. \nonumber \\
& + & \int_{0}^{t}||S_{n}(r) - S_{\infty}(r)||_{\infty}||M(r)^{\tra}\Phi(r)M(r)\Psi(r)||_{\infty}dr \nonumber \\
& + & \int_{0}^{t}||\Phi(r) - \Sigma_{n}(r)^{-1}||_{\infty}||S_{n}(r)M(r)^{\tra}M(r)\Psi(r)||_{\infty}dr
\end{eqnarray}
To proceed further we need some additional estimates.
Using (\ref{obs}) and writing $J(t): = M(t)^{\tra}\Phi(t)G(t)$ for all $0 \leq t \leq T$, we find that
\bean & & \E\left(\left| \int_{0}^{t}(S_{n}(r) - S_{\infty}(r))J(r)dZ(r)\right|\right)\\
& \leq & \E\left(\left| \int_{0}^{t}(S_{n}(r) - S_{\infty}(r))J(r)C(r)Y(r)dr \right|\right)\\
& + & \E\left(\left| \int_{0}^{t}(S_{n}(r) - S_{\infty}(r))J(r)D(r)dW_{a}(r)\right|\right)\\
& + & \E\left(\left| \int_{0}^{t}\int_{\R}(S_{n}(r) - S_{\infty}(r))J(r)D(r)\tilde{N}(dy,dr)\right|\right)\\
& \leq & \int_{0}^{t}||S_{n}(r) - S_{\infty}(r)||_{\infty}||J(r)||_{\infty}||C(r)||_{\infty}||Y(r)||_{1}dr\\
& + & \left(\E\left(\left| \int_{0}^{t}(S_{n}(r) - S_{\infty}(r))J(r)D(r)dW_{a}(r)\right|^{2}\right)\right)^{\frac{1}{2}}\\
& + & 2\int_{\R}|y|\nu(dy)\int_{0}^{t}||S_{n}(r) - S_{\infty}(r)||_{\infty}||J(r)||_{\infty}||C(r)||_{\infty}dr,\eean
and by It\^{o}'s isometry,
\bean & & \E\left(\left| \int_{0}^{t}(S_{n}(r) - S_{\infty}(r))J(r)D(r)dW_{a}(r)\right|^{2}\right)\\
& = & \int_{0}^{t}||(S_{n}(r) - S_{\infty}(r))J(r)D(r)a^{\frac{1}{2}}||_{\mbox{hs}}^{2}dr\\
& \leq & \tr(a)\int_{0}^{t}||(S_{n}(r) - S_{\infty}(r)||^{2}_{\infty}||J(r)||_{\infty}^{2}||D(r)||_{\infty}^{2}dr. \eean
The other term that involves $\E$ is dealt with in the same manner.
We incorporate these estimates into (\ref{last}) and also make use of (\ref{est}) to conclude that there exist bounded measurable positive functions $t \rightarrow H_{j}(t)$ defined on $[0,T]$ for $j=1,\ldots ,6$ so that

\bean & &  ||\Psi(t) - \widehat{Y_{n}}(t)||_{1}\\ & \leq & \int_{0}^{t} H_{1}(r)||\Psi(r) - \widehat{Y_{n}}(r))||_{1}dr
 +  \int_{0}^{t}H_{2}(r)||S_{n}(r) - S_{\infty}(r)||_{\infty}dr \\
& + & \int_{0}^{t}H_{3}(r)||\Phi(r) - \Sigma_{n}(r)^{-1}||_{\infty}dr
 +  \left( \int_{0}^{t}H_{4}(r)||S_{n}(r) - S_{\infty}(r)||_{\infty}dr\right)^{\frac{1}{2}}\\
& + & \left( \int_{0}^{t}H_{5}(r)||\Phi(r) - \Sigma_{n}(r)^{-1}||_{\infty}dr\right)^{\frac{1}{2}}
 +  \int_{|y| > n}|y|\nu(dy)\int_{0}^{t}H_{6}(r)dr. \eean

Then by Gronwall's lemma
\bean ||\Psi(t) - \widehat{Y_{n}}(t)||_{1} & \leq & \exp{\left(\int_{0}^{t}H_{1}(r)dr\right)}\left[\int_{0}^{t}H_{2}(r)||S_{n}(r) - S_{\infty}(r)||_{\infty}dr\right. \\
& + & \int_{0}^{t}H_{3}(r)||\Phi(r) - \Sigma_{n}(r)^{-1}||_{\infty}dr \\
& + & \left( \int_{0}^{t}H_{4}(r)||S_{n}(r) - S_{\infty}(r)||_{\infty}dr\right)^{\frac{1}{2}}\\
& + & \left( \int_{0}^{t}H_{5}(r)||\Phi(r) - \Sigma_{n}(r)^{-1}||_{\infty}dr\right)^{\frac{1}{2}}\\
& + &  \left. \int_{|y| > n}|y|\nu(dy)\int_{0}^{t}H_{6}(r)dr \right], \eean
and the result follows by using dominated convergence (which is justified by using (\ref{est})) and the fact that $\lim_{n \rightarrow \infty} \int_{|y| > n}|y|\nu(dy) =  0$. $\hfill \Box$

\vspace{5pt}

When the observation noise has fully infinite variance (\ref{KBFinf}) collapses to a linear first order differential equation:

\begin{equation} \label{KBFinf1}
\frac{d\widehat{Y}(t)}{dt}  =  A(t)\widehat{Y}(t).
\end{equation}

\section{Numerical Results}

In this section we look at three different methods for filtering an infinite variance observation process in a linear framework.  We work in the one-dimensional case using (\ref{KBFinf1}) and compare our results with those of \cite{AF} and \cite{BM}.

The example we will look at will be that of infinite variance $\alpha$-stable noisy obervations of a mean reverting Brownian motion, i.e.
\bean
Y(t) &=&Y_0  -\int_0^t Y(t)dt+ X_1(t)\nonumber \\
dZ(t) &=& Y(t)dt + dX_2(t) \nonumber
\eean
where $(X_1(t),t\geq 0)$ is a standard Brownian motion, $Y_0$ is a Gaussian random variable and $(X_2(t),t\geq 0)$ is a symmetric $\alpha$-stable process with dispersion parameter $1$ and $\alpha \in (1,2)$, so that $\E(e^{iuX_{2}(t)}) = e^{-t|u|^{\alpha}}$ for each $u \in \R, t \geq 0$.  We set the signal to be $0$ at time $0$ and assume the observations occur at a fixed rate with interarrival time $0.01$ and expiry $T=10$. Finally we take $100,000$ Monte-Carlo simulations to estimate the  mean square error of the filter from the system, after which a median value of these were taken, the results are summarised in Table \ref{table:lin}.

\begin{table}[h]
\caption{A Comparison of Filters}
\centering
\begin{tabular}{c c c c}
\hline\hline
$\alpha$ & No Filter Error & LBM Filter Error & AF Filter Error\\ [0.5ex]
\hline
1.1 & 0.6601 & 0.7819 & 0.7445\\
1.5 & 0.6593 &  0.6819 & 0.7507\\
1.9 & 0.6603 & 0.6446  & 0.7508\\ [1ex]
\hline
\end{tabular}
\label{table:lin}
\end{table}

Here ``No filter'' refers to the filter of this paper and LBM and AF are those developed in \cite{BM} and \cite{AF}, respectively.

The following figure provides a graphical representation of one instance of the above calculations. It shows that the filter of Le Breton and Musiela is quite prone to large jumps in the observations, causing a higher error in this instance. The filter of Ahn and Feldman remains  relatively close to $0$, but having no filter provides the lowest error. Whilst this may seem counter-intuitive it should be noted that when $\alpha = 1.1$ the noise in the observation process is very wild, and so the observation process bears almost no resemblance to the system process.

\begin{figure}[h]
\centering
\includegraphics[scale=0.30]{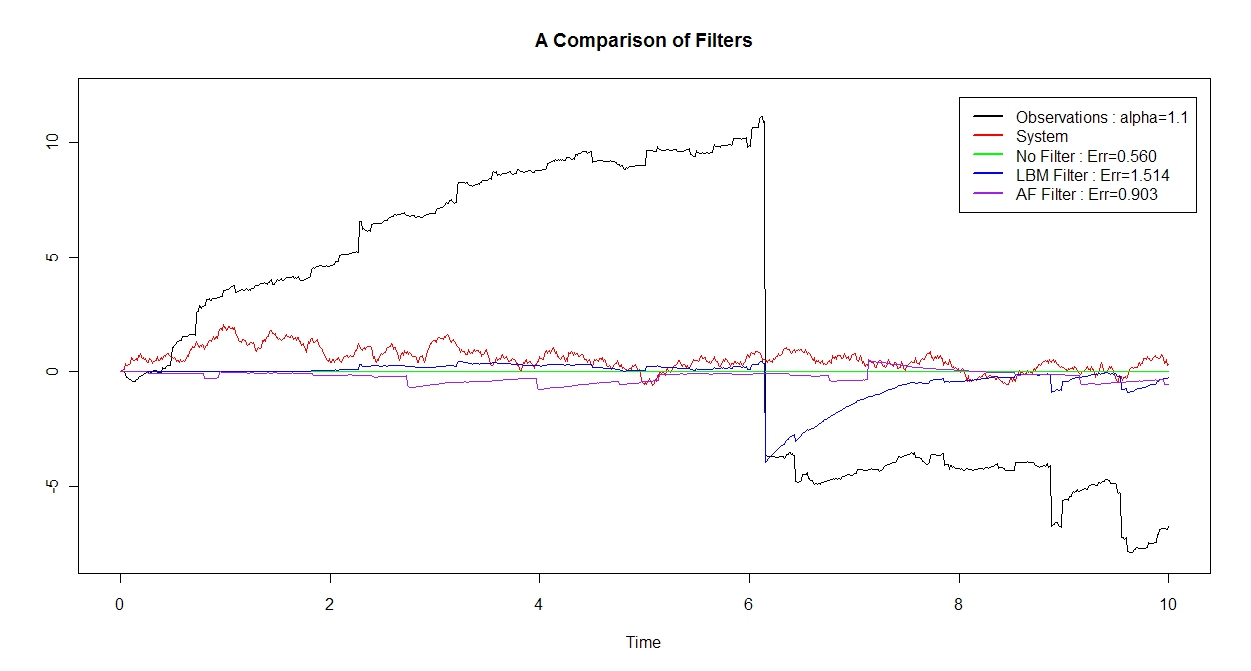}
\end{figure}

\pagebreak

{\it Acknowlegements.} We thank John Biggins and Dan Crisan for very helpful comments.


\begin{thebibliography}{99}

\bibitem{AF} H.Ahn, R.E.Feldman, Optimal filtering of a Gaussian signal in the prescence of L\'{e}vy noise, {\it Siam J. Appl. Math.} {\bf 60}, 359-69 (1999)

\bibitem{App} D.Applebaum, {\it L\'{e}vy Processes and Stochastic Calculus} (second edition), Cambridge University Press (2009)

\bibitem{App1} D.Applebaum, Extending stochastic resonance for neuron models to general Levy noise, {\it IEEE Trans. Neural Netw.}  {\bf  20} 1993-6 (2009)

\bibitem{BC} A.Bain, D.Crisan, {\it Fundamentals of Stochastic Filtering}, Stochastic Modelling and Applied Probability {\bf 60}, Springer (2009)

\bibitem{Bar} F.Bartumeus, L\'{e}vy processes in animal movement: an evolutionary hypothesis, {\it Fractals} {\bf 15}, 1-12 (2007)

\bibitem{Blac} S.Blackwood, {\it Stochastic Filtering with L\'{e}vy Processes}, University of Sheffield PhD thesis (2014)

\bibitem{BM} A. Le Breton, M. Musiela, A generalisation of the Kalman filter to models with infinite variance, {\it Stoch. Proc. App.} {\bf 47}, 75-94 (1993)

\bibitem{Cip} B.Cipra, Engineers look to Kalman filtering for guidance, {\it SIAM News} {\bf 26}, No.5 (1993).

\bibitem{CoTa} R.Cont, P.Tankov, {\it Financial Modelling with Jump Processes}, Chapman and Hall/CRC (2004)

\bibitem{Dav} M.H.A.Davis, {\it Linear Estimation and Stochastic Control}, Chapman and Hall Mathematical Series, Chapman and Hall, London (1977)


\bibitem{Halm} P.R.Halmos, {\it A Hilbert Space Problem Book} (second edition), Springer-Verlag New York Inc. (1982)

\bibitem{KB} R.E.Kalman, R.S.Bucy, New results in linear filtering and prediction theory, {\it Trans. ASME Ser. D. J. Basic Engrg.} {\bf 83}, 95�108 (1961)

\bibitem{KB} M.L.Kleptsyna, A.Le Breton, Extension of the Kalman-Bucy filter to elementary linear systems with fractional Brownian noises, {\it Statist. Infer. Stoch. Proc.} {\bf 5} 249-71 (2002)

\bibitem{MBP} T.Meyer-Brandis, F.Proske, Explicit solution of a non-linear filtering problem for L\'{e}vy processes with application to finance, {\it Appl. Math. Optim.} {\bf 50}, 119-34 (2004)

\bibitem{Ox} B.{\O}ksendal, {\it Stochastic Differential Equations} (sixth edition), Springer-Verlag (2003).

\bibitem{KP} A.Patel, B.Kosko, Stochastic resonance in continuous and spiking neuron models with Levy
noise, {\it IEEE Trans. Neural Netw.}, {\bf 19} 1993-2008, (2008)

\bibitem{PS} S.Popa, S.S.Sritharan, Nonlinear filtering of It\^{o}-L\'{e}vy stochastic differential equations with continuous obervations, {\it Commun. Stoch. Anal.} {\bf 3}, 313-30 (2009)

\bibitem{SI} D.Sornette, K.Ide, The Kalman-L\'{e}vy filter, {\it Physica D} {\bf 151}, 142-74 (2001)

\bibitem{Xi} J.Xiong, {\it An Introduction to Stochastic Filtering Theory}, Oxford Graduate Texts in Mathematics {\bf 18}, Oxford University Press (2008)


\end{thebibliography}
\end{document}